\input amstex
\documentstyle{amsppt}
\magnification 1200
\vcorrection{-9mm}
\input epsf

\def\R{\Bbb R}
\def\C{\Bbb C}
\def\CP{\Bbb CP}
\def\bS{\Bbb S}

\def\cQ{\Cal Q}
\def\cSB{\Cal{SB}}
\def\cB{\Cal B}
\def\bb{\bold b}
\def\Int{\operatorname{Int}}
\def\sgn{\operatorname{sgn}}
\def\lk{\operatorname{lk}}

\topmatter
\title
       $\Bbb C$-boundary links up to six crossings
\endtitle
\author
       S.~Yu.~Orevkov
\endauthor
\abstract
An oriented link is called {\it $\C$-boundary} if it is realizable as
$(\partial B,A\cap\partial B)$ where $A$ is an algebraic curve in $\C^2$ and $B$
is an embedded $4$-ball. This notion was introduced by Michel Boileau and Lee Rudolph
in 1995.
In a recent joint paper with N.\,G.~Kruzhilin we gave a complete classification of
$\C$-boundaries with at most 5 crossings. In the present paper a more regular method
of construction of $\C$-boundaries is proposed and the classification is extended
up to 6 crossings.
\endabstract
\address
        Steklov Mathematical Institute of Russian Academy of Sciences, Moscow, Russia.
\endaddress

\address
        IMT, l'universit\'e Paul Sabatier, 118 route de Narbonne, Toulouse, France.
\endaddress
\endtopmatter

\def\sectConstr  {2}
\def\sectSixCr   {3}
\def\sectIY      {3.1}
\def\sectZapret  {3.2}
\def\sectTab     {3.3}
\def\sectConclude{4}
\def\sectSqueezed   {4.1}
\def\sectSliceKnots {4.2}
\def\sectQuestions  {4.3}
 
\def\propConstr    {\sectConstr.1}
\def \corConstr    {\sectConstr.2}
\def\exaQQ    {\sectConstr.3}
\def\exaQQii  {\sectConstr.4}
\def\exaQQiii {\sectConstr.5}
\def\propConstrStrong  {\sectConstr.6}
\def\exaFive  {\sectConstr.7}

\def\thKM    {\sectSixCr.1}
\def\propB   {\sectSixCr.2}
\def\propC   {\sectSixCr.3}
\def\propD   {\sectSixCr.4}

\def\propSqueezed {\sectConclude.1}
\def\queConstr  {\sectConclude.2}
\def\queSliceI  {\sectConclude.3} 
\def\queSliceII {\sectConclude.4} 
\def\queStrong  {\sectConclude.5} 
\def\queCancel  {\sectConclude.6} 
\def\queSplit   {\sectConclude.7} 

\def\figBA      {1}
\def\figIsot    {2}
\def\figIsotII  {3}
\def\figIsotIII {4}
\def\figCor     {5}
\def\figQQ      {6}
\def\figFive    {7}

\def\figD     {8}
\def\figDD    {9}
\def\figMain {10}

\def\figSliceKnots {11}

\def\tabA   {1}
\def\tabB   {2}
\def\tabC   {3}

\def\refBO  {1}
\def\refBR  {2}
\def\refFLL {3}
\def\refFS  {4}
\def\refFW  {5}
\def\refKM  {6}
\def\refKO  {7}
\def\refLMknot {8}
\def\refLMlink {9}
\def\refM      {10}
\def\refMo     {11}
\def\refO      {12}
\def\refRuTop  {13}
\def\refRuFancy {14}
\def\refRuObstr {15}

\document

\head 1. Introduction
\endhead

Boileau and Rudolph [\refBR]
defined {\it $\C$-boundary} as an oriented link in
the $3$-sphere $\bS^3$ which can be realized as $(\partial B,A\cap\partial B)$
where $A$ is an algebraic curve in $\C^2$ and $B$ is an embedded $4$-ball
(with the boundary orientation of $A\cap\partial B$ induced from $A\cap B$).
It is observed in [\refBR] that Kronheimer-Mrowka Theorem [\refKM] (former
Thom Conjecture) implies that some knots are not $\C$-boundaries, for example,
the figure-eight knot $4_1$.

If $B$ is pseudo-convex (for example, a standard $4$-ball), then
$(\partial B,A\cap\partial B)$ is a {quasipositive link} by [\refBO].
In [\refKO] we observed that there exist non-quasipositive $\C$-boundaries,
the simplest ones being $L\,\#\,{-L^*}$ where $L$ is a quasipositive link.
Here and below, $L^*$ denotes the mirror image of $L$, and
$-L$ denotes $L$ with the opposite orientation. In [\refKO] we
also corrected some errors in [\refBR] and observed that
it is often more efficient to apply the so-called
Immersed Thom Conjecture in order to obtain restrictions on the $\C$-boundaries
(see details in [\refKO, \S3] and in \S3 here).

It is natural to distinguish the case when $L$ is realizable as $A\cap\partial B$
as above, and moreover $A\setminus B$ is connected.
Such links are called in [\refKO] {\it strong $\C$-boundaries}.
Examples of non-strong $\C$-boundaries are given in [\refKO]. The simplest one
is the split sum $2_1\sqcup 2_1^*$ (we denote the positive
two-component Hopf link by $2_1$).

We gave in [\refKO] a complete list of $\C$-boundaries and strong
$\C$-boundaries with at most 5 crossings (including split and composite links).
In the present paper a complete list of $\C$-boundaries with 6 crossings is given
(see \S\sectSixCr). For two of them it remains unknown if they are strong
$\C$-boundaries. Very elementary constructions were used in [\refKO].
They are not enough for links with 6 crossings, and in \S\sectConstr\ we present
a more systematic way to construct (strong) $\C$-boundaries. In contrary,
the fact that some links are not (strong) $\C$-boundaries, is proven here by the same
methods as in [\refKO]. In \S\sectSliceKnots\ we show that 10 (among 15)
slice knots with $\le 9$ crossings are $\C$-boundaries.
Some open questions are posed in \S\sectQuestions. 

\smallskip
{\bf Acknowledgement.} I thank Michel Boileau and the anonymous referee for valuable remarks.


\head\sectConstr. A construction of $\C$-boundaries
\endhead

\subhead\sectConstr.1. Notation and terminology
\endsubhead
We denote the segment $[0,1]$ by $I$ and we set $I^2=I\times I$.
Let $L$ be an oriented link.
A {\it band attached to} $L$ is an embedding $b:I^2\to\bS^3$ such that
$b|_{I\times\partial I}\to L$ is an orientation preserving embedding
where the orientation of $I\times\partial I$ is inherited from the boundary
orientation of $\partial(I^2)$
(see Figure~\figBA). Sometimes, abusing the language, we use this term
for $b(I^2)$ rather than for $b$. The result of the {\it band-move}
(or {\it band surgery}) on $L$ along $b$ is
$$
    L_b=\big(L\setminus b(I\times\partial I)\big)\cup b(\partial I\times I)
$$
with the orientation coherent with $L$ (see Figure~\figBA). We say that two
bands $b,b'$ are disjoint if $b(I^2)\cap b'(I^2)=\varnothing$.
If $\bb$ is a collection of pairwise disjoint bands attached to $L$, let $L_\bb$
denote the result of simultaneous band-moves on $L$ along all these bands.

\if01{
A band $b$ attached to a link $L$ is called {\it trivial} if $L_b$ is a split
sum $L'\sqcup O$ where $O$ is an unknot such that $O\not\in\cap L\ne\varnothing$
(see the right panel in Figure~\figBA). In this case, $O\cap L$ is called
the {\it vanishing arc} of $(L,\bb)$.
A collection $\bb=\{b_1,\dots,b_n\}$ of pairwise disjoint bands attached to $L$ is
{\it trivial} if $L_\bb$ is a split sum $L'\sqcup O_1\sqcup\dots\sqcup O_n$ where
each $O_i$ is an unknot composed of an arc of $\partial b_i$ and an arc of $L$.
When such $\bb$ is a subset of a larger collection $\bb'$ of pairwise disjoint bands,
we say that $\bb$ is a {\it trivial subset of} $\bb'$ if $L_\bb$ admits the above
splitting such that all the $O_i$ are disjoint from the bands from $\bb'\setminus\bb$.
}\fi

A band $b$ attached to a link $L$ is called {\it trivial} if $L_b$ is a split
sum $L'\sqcup O$ where $O$ is an unknot composed of an arc of
$\partial b$ and an arc of $L$ (see the right panel in Figure~\figBA).
In this case, $O\cap L$ is called the {\it vanishing arc} of $(L,\bb)$.
A collection $\bb=\{b_1,\dots,b_n\}$ of pairwise disjoint bands attached to $L$ is
{\it trivial} if all $b_i$ are trivial for $L$ and their vanishing arcs are
pairwise disjoint.
When such $\bb$ is a subset of a larger collection $\bb'$ of pairwise disjoint bands,
we say that $\bb$ is a {\it trivial subset of} $\bb'$, if it is a trivial collection
of bands and all its vanishing arcs are disjoint from the bands belonging
to $\bb'\setminus\bb$.

\midinsert
\centerline{\epsfxsize=100mm \epsfbox{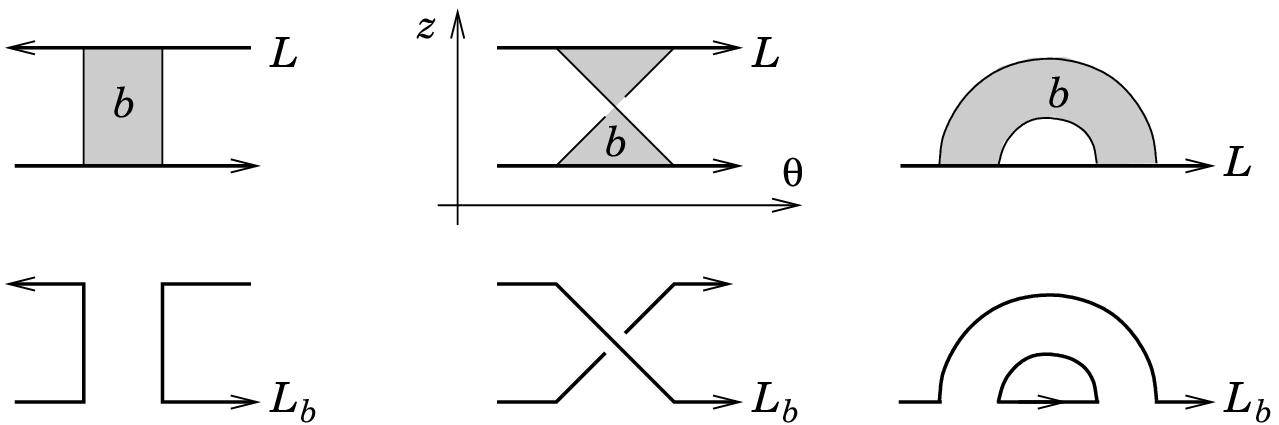}}
\botcaption{Figure \figBA} A band $b$ attached to a link $L$. The link $L_b$.
         \hbox to 20mm{}
         In the middle: positive $\theta$-monotone band. On the right: trivial band.
\endcaption
\endinsert

Now we identify $\bS^3$ with a one-point compactification of $\R^3$ with standard
coordinates $(x,y,z)$.
Let us fix cylindrical coordinates $(r,\theta,z)$ in $\R^3$, that is
$x=r\cos\theta$, $y=r\sin\theta$. We say that
a link $L$ is {\it $\theta$-monotone} 
if $r|_L>0$ and $d\theta$ is positive on $L$.
We also call such links {\it geometric closed braids}.

Let $L$ be a $\theta$-monotone link. A band $b$ attached to $L$ is called
{\it positive $\theta$-monotone} if $L_b$ is $\theta$-monotone and the
projection of some neighborhood of $L\cup b(I^2)$ to the $(\theta,z)$-plane looks as
in the middle panel of Figure~\figBA.

An isotopy $\{L_t\}_{t\in I}$ is called {\it $\theta$-monotone} if each link $L_t$
is $\theta$-monotone.


\subhead\sectConstr.2. The construction
\endsubhead
Let $L$ be a $\theta$-monotone link which is the braid closure of a quasipositive
(maybe, trivial) braid, and let $\bb$ be a collection of pairwise disjoint
bands attached to $L$.
Suppose that there exists a sequence of links $L=L_0,L_1,\dots,L_n=L'$
and a sequence of collections of bands $\bb=\bb_0,\bb_1,\dots,\bb_n=\bb'$
($\bb_i$ is attached to $L_i$) such that for each $i=1,\dots,n$, one of
the following two conditions holds:
\roster
\item
      $L_i\cup\bb_i$ is obtained from $L_{i-1}\cup\bb_{i-1}$ by an isotopy
      which restricts to a $\theta$-monotone isotopy between $L_{i-1}$ and $L_i$ or
\item
      $L_i=(L_{i-1})_{c_i}$ and $\bb_i=\bb_{i-1}$
      where $c_i$ is a positive $\theta$-monotone band attached to $L_{i-1}$
      (see Figure~\figBA) disjoint from $\bb_{i-1}$.
\endroster


\proclaim{ Proposition \propConstr }
If $\bb'$ is trivial, 
then $L_\bb$ is a $\C$-boundary.
\endproclaim

\demo{ Proof }
%
By construction, $L'$ is the braid closure of a quasipositive braid.
Hence, by Rudolph's theorem [\refRuTop], $(\bS^3,L')$ is diffeomorphic
to $(\partial B',A\cap\partial B')$, $B'=\Delta'\times\Delta''$,
where $\Delta'$ and $\Delta''$
are some disks in $\C$, and $A$ is a complex algebraic curve in $\C^2$ disjoint
from $\Delta'\times\partial\Delta''$.
Moreover, the arguments in [\refRuTop] imply that there exist embedded
disks $\Delta=\Delta_0\subset\Delta_1\subset\dots\subset\Delta_n=\Delta'$
(with $\Delta_{i-1}\subset\Int\Delta_i$) such that
$(B_i,A\cap\partial B_i)\approx(\bS^3,L_i)$, $i=0,\dots,n$, where
$B_i=\Delta_i\times\Delta''$, and each time when $L_i$ is obtained from
$L_{i-1}$ by a positive $\theta$-monotone band-move, 
the annulus $\Delta_i\setminus\Delta_{i-1}$ contains a ramification point
of $\pi|_{A\cap B'}$ where $\pi:B'=\Delta'\times\Delta''\to\Delta'$ is
the projection onto the first factor.

We include the disk boundaries $\partial\Delta_0,\dots,\partial\Delta_n$
into a continuous family. Namely, consider a diffeomorphism
$f:\bS^1\times I\to\Delta'\setminus\Int\Delta$, such that
$\partial\Delta_i=f_{t_i}(\bS^1)$ for some numbers $t_i$
in the range $0=t_0<t_1<\dots<t_n=1$ (here $f_t:\bS^1\to\Delta'$ is defined by
$f_t(p)=f(p,t)$).

For each band $b$ from $\bb$, the isotopies between $\bb_{i-1}$ and $\bb_i$
(which are identical when $\bb_i=(\bb_{i-1})_{c_i}$) can be glued together
into a single isotopy 
$\{b_t:I^2\to B'\}_{t\in I}$
such that $\pi(b_t(I^2))\subset f_t(\bS^1)$.
Let $\{h_t:I^3\to B'\}_{t\in I}$ be an isotopy obtained from $\{b_t\}$
by replacing each band $b_t:I^2\to\pi^{-1}(f_t(\bS^1))$ with it thickening
$h_t:I^3\to\pi^{-1}(f_t(\bS^1))$ as shown in Figure~\figIsot.
Finally, we modify this isotopy near $t=1$ as shown in Figure~\figIsotII\
(still keeping the notation $h_t$ for it).

\midinsert
\centerline{\epsfxsize=55mm \epsfbox{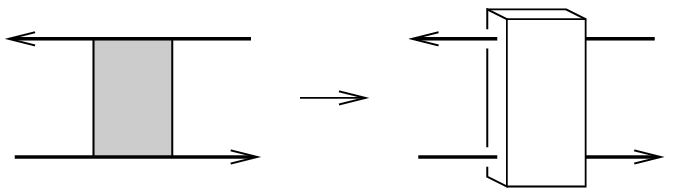}}
\botcaption{Figure \figIsot} A thickened band.
\endcaption
\endinsert

\midinsert
\centerline{\epsfxsize=90mm \epsfbox{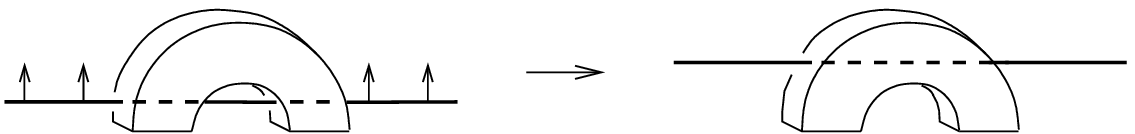}}
\botcaption{Figure \figIsotII} The final isotopy of a thickened trivial band.
\endcaption
\endinsert

We obtain an embedding of the 4-cube 
$h=h^{(b)}:I^4=I^3\times I\to B'$, $h(p,t)=h_t(p)$.
The intersection of $h(I^4)$ with $B=\Delta\times\Delta''$
is the image of one of the eight faces of $I^4$,
that is $h(I^4)\cap B=h(I^3\times 0)$. Thus the following set $\hat B$
(appropriately smoothed) is an embedded 4-ball:
$$
   \hat B = B\cup \bigcup_{b\in\bb} h^{(b)}(I^4).
$$
Let us show that
$(\partial\hat B,A\cap\partial\hat B)\approx(\bS^3,L_\bb)$.
Indeed, when the 4-cube $h^{(b)}(I^4)$ is attached to $B$,
the following surgery is applied to the boundary link:
the 3-ball $h^{(b)}(I^4\times 0)$, which is a face of the attached 4-ball,
is replaced by the union of the seven other faces. One can see in Figure~\figIsotIII\
that this operation coincides with the band-move on $L$ along $b$
(cf.~Figure~\figBA).
\qed\enddemo

\midinsert
\centerline{\epsfxsize=90mm \epsfbox{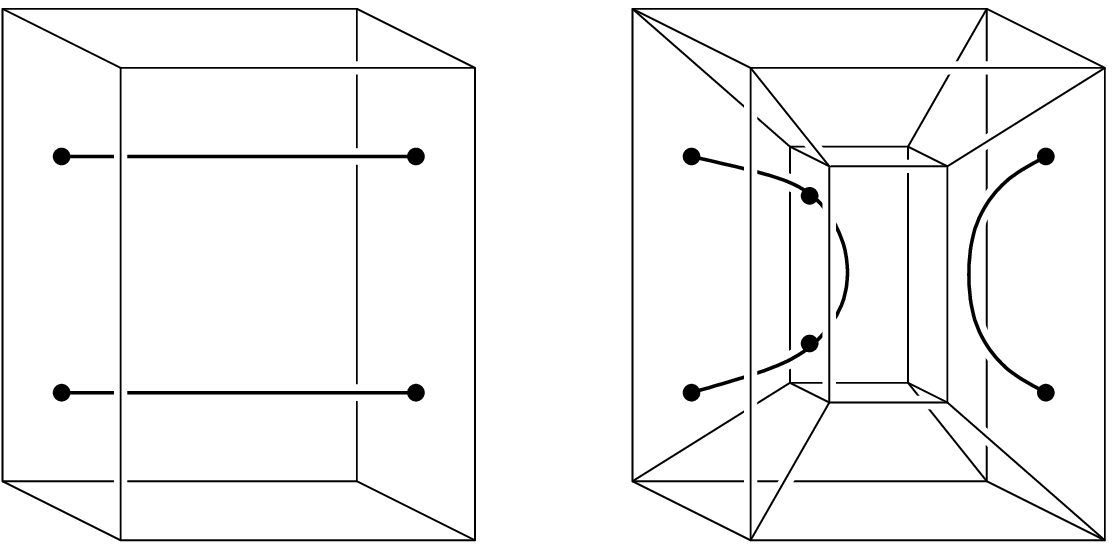}}
\botcaption{Figure \figIsotIII} Replacement of one face of $I^4$ by seven other faces.
\endcaption
\endinsert

\proclaim{ Corollary \corConstr }
If $\bb'=\bb'_0\cup\bb'_1$ where $\bb'_0$ is a trivial subset of $\bb'$
and each band in $\bb'_1$ is positive $\theta$-monotone,
then $L_\bb$ is a $\C$-boundary.
\endproclaim

\demo{ Proof }
We continue the sequences $(L_1,\dots\;)$ and $(\bb_1,\dots\;)$ as shown in
Figure~\figCor. Then each positive $\theta$-monotone band is transformed into
a trivial band.
\qed\enddemo

\midinsert
\centerline{\epsfxsize=100mm \epsfbox{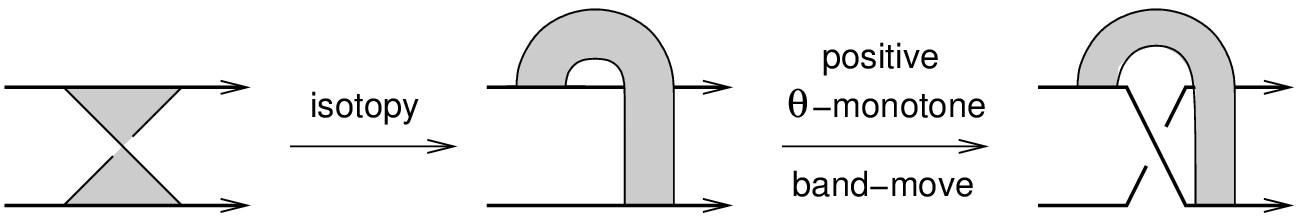}}
\botcaption{Figure \figCor} From a positive $\theta$-monotone band to a trivial band.
\endcaption
\endinsert

\if01{
If $\bb$ is trivial, then $L'_{\bb'}$ is a $\C$-boundary.

\demo{ Proof }
The proof is almost the same as for Proposition~\propConstr, but the
embedded 4-cubes $h^(b)(I^4)$ are attached from the interior side to
$\partial B'$, and the resulting link $(\bS,L'_{\bb})$ is realized by
$(\partial\hat B',A\cap\partial\hat B')$ where
$\hat B' = B'\setminus \bigcup_{b\in\bb'} h^{(b)}(I^4)$.
\qed\enddemo

The both constructions can be evidently combined as follows.

Let $\bb$ be a collection of pairwise disjoint bands attached to
a $\theta$-monotone link $L$. Suppose that $\bb=\bb_1\cup\bb_2$,
$\bb_1\cap\bb_2=\varnothing$, and there exist quasipositive
cobordisms from $L''\cup\bb''$ to $L\cup\bb_1$ and from $L\cup\bb_2$
to $L'\cup\bb'$ where $L''$ is the braid closure of a quasipositive braid
and both $\bb'$ and $\bb''$ are trivial.
Then $L_\bb$ is a $\C$-boundary.
}\fi

\medskip\noindent
{\bf Example \exaQQ.}
It is almost evident (see [\refKO, \S2]) that $L\sqcup-L^*$ and $L\,\#\,{-L^*}$ are
$\C$-boundaries, if $L$ is a quasipositive link (in 
$L\,\#\,{-L^*}$, some component of $L$ is supposed to be joint with its mirror image).
This fact can be also obtained as a simplest application of Proposition~\propConstr.
Indeed, let $\beta$ be a quasipositive braid whose braid closure is $L$.
We start with the geometric braid with bands shown in Figure~\figQQ,
and apply positive $\theta$-monotone band-moves to the $\beta^{-1}$
half so that it is transformed into a trivial braid.

\midinsert
\centerline{\epsfxsize=80mm \epsfbox{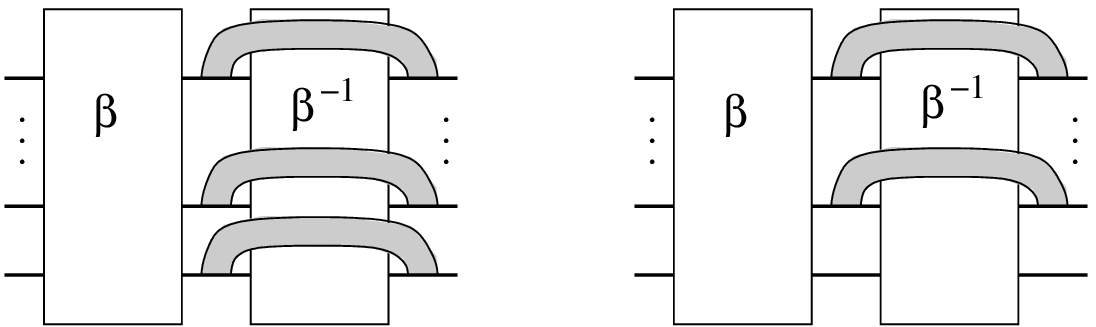}}
\botcaption{Figure \figQQ} $\C$-boundary realization of 
    $L\,\sqcup-L^*$ (on the left) and $L\,\#\,{-L^*}$ (on the right)
     based on Proposition~\propConstr.
\endcaption
\endinsert

\medskip\noindent
{\bf Example \exaQQii.}
Similarly to Example~\exaQQ, we obtain a $\C$-boundary realization of
$L\,\sqcup-L_1^*$ and $L\,\#\,{-L_1^*}$ where $L$ and $L_1$ are the braid closures
of $\beta=\beta_1\beta_2$ and $\beta_1$ respectively for quasipositive braids
$\beta_1$ and $\beta_2$: just replace $\beta^{-1}$ by $\beta_1^{-1}$ in Figure~\figQQ.

\medskip\noindent
{\bf Example \exaQQiii.}
Similarly to the previous examples, we obtain a $\C$-boundary realization of
$L\,\#_{(K,-K^*)}\,{-L^*}$ where $L$ is the braid closure of an $n$-braid $\beta$
belonging to the monoid $M$ generated by
$\{\sigma_1,\dots,\sigma_{k-1},\sigma_k^{-1},\dots,\sigma_{n-1}^{-1}\}$
for some $k\le n$, and $K$ is the component of $L$ containing the $k$-th
string of $\beta$. To this end, we start with a configuration as in Figure~\figQQ\
but with bands attached to all the strings except the $k$-th one.
Then, by a sequence of positive $\theta$-monotone band-moves we eliminate
all the negative crossings and obtain a trivial collection of bands attached to
a positive braid.
In particular $4_1\# 4_1$ is a $\C$-boundary: $4_1$ is the
braid closure of $(\sigma_1\sigma_2^{-1})^2$.

The same construction (and its generalization in the spirit of Example~\exaQQii)
also applies if $M$ is replaced by the monoid generated by
$$
   \{a\sigma_1a^{-1}\mid a\in \iota_0(B_k)\} \cup
   \{a\sigma_k^{-1}a^{-1}\mid a\in \iota_{k-1}(B_{n-k+1})\}
$$
where $\iota_p:B_m\to B_n$ (for $m+p\le n$) is the homomorhism of the braid groups
which takes $\sigma_i$ of $B_m$ to $\sigma_{i+p}$ of $B_n$.


\subhead\nofrills\sectConstr.3. When the obtained $\C$-boundary is strong?
\endsubhead\;\;
Let the setting be as in \S\sectConstr.2.
Suppose that $\bb'$ is trivial, and hence $L_\bb$ is a $\C$-boundary
by Proposition~\propConstr. Let us explain how to check if it is strong or not.


Let $G$ be the graph whose vertices correspond to all components
(which are open intervals or circles) of $L_i\setminus\bb_i$ for all $i=0,\dots,n$,
and the edges are defined as follows.
If $L_i\cup\bb_i$ is
obtained from $L_{i-1}\cup\bb_{i-1}$ by a $\theta$-monotone isotopy, then
the corresponding pairs of vertices are connected by an edge of $G$.
If $L_i=(L_{i-1})_{c_i}$ and $\bb_i=\bb_{i-1}$,
then all the four vertices corresponding the components of
$L_{i-1}\setminus\bb_{i-1}$ and
$L_i\setminus\bb_i$ adjacent to $c_i$ are connected to each other by edges of $G$.

\proclaim{ Proposition \propConstrStrong }
If each vertex of $G$ corresponding to a vanishing arc of $(L',\bb')$
is connected by a path with a vertex corresponding to a non-vanishing arc of
$(L',\bb')$, then the link $L_\bb$ is a strong $\C$-boundary.
\endproclaim

Note that the hypothesis of Proposition~\propConstrStrong\ is equivalent
to the fact that the sequences 
 $L_1,\dots,L_n$ and $\bb_1,\dots,\bb_n$
can be extended so that
the graph $G$ becomes connected.
One can check that the hypothesis of Proposition~\propConstrStrong\
does not hold for $L\sqcup-L^*$ in Example~\exaQQ\
(cf.~Question~\queStrong\ below).

\medskip\noindent
{\bf Example \exaFive.} (A new $\C$-boundary realization of $5_1^2$.)
In [\refKO, Prop.~6.3] we proved that the link $5_1^2$ (the mirror image of
L5a1 in [\refLMlink]) is a strong $\C$-boundary.
Corollary~\corConstr\ provides another proof which seems to be simpler
(see Figure~\figFive).

\midinsert
\centerline{\epsfxsize=110mm \epsfbox{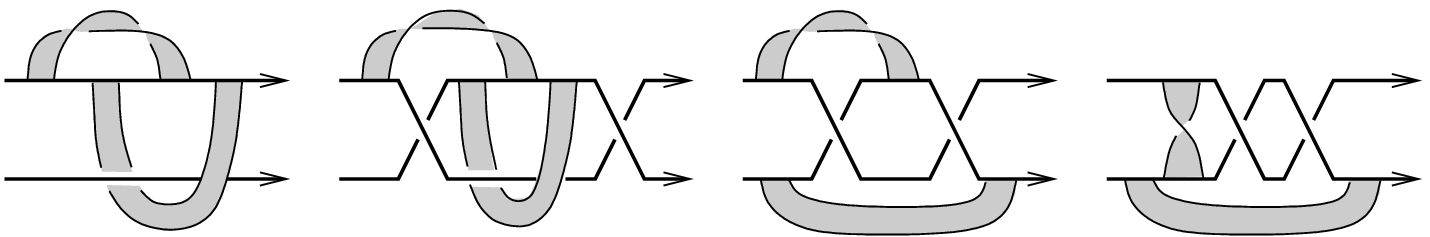}}
\botcaption{Figure \figFive} A new strong $\C$-boundary realization of $5_1^2$.
\endcaption
\endinsert


\head\sectSixCr. Links with six crossings
\endhead

As in [\refKO], we denote the class of quasipositive links,
the class of $\C$-boundaries, and the class of strong $\C$-boundaries by $\cQ$,
$\cB$, and $\cSB$ respectively.
In this section, for each link with crossing number 6,
we determine if it belongs to these classes (with two exceptions for
the class $\cSB$).
In Tables~\tabA--\tabC\ we give the answers for all such links which
do not have an unknot split component,
however, the answers for links of the form
$L\sqcup O\sqcup\dots\sqcup O$ with $6$ crossings easily follow (see [\refKO, \S6]).
More explanations about Tables~\tabA--\tabC\ are given in \S\sectTab.

\if01{
The implication $L\in\Cal C\Rightarrow L\sqcup O\in\Cal C$
($\Cal C$ is $\cQ$, $\cSB$, or $\cB$) is evident.
The reverse implication is known only for $\cQ$ (see [\refO]), but
the nature of our proofs is such that each time when we prove that
$L\not\in\Cal C$ (where $\Cal C$ is the class $\cSB$ or $\cB$)
the same arguments can be easily
adapted to prove that $L\sqcup O\sqcup\dots\sqcup O\not\in\Cal C$.
}\fi

\midinsert
\hbox to 100mm{\hfill Table \tabA. Prime links.}
\medskip
\centerline{
\vbox{\offinterlineskip
\def\h {height2pt&\omit&&\omit&&\omit&&\omit&&\omit&&\omit&&\omit&&\omit&&\omit&\cr}
    
\def\1{\bar1}\def\2{\bar2}\def\3{\bar3}\def\4{\bar4}\def\*{${}^*$}

\hrule
\halign{&\vrule#&\strut\;\hfil#\hfil\,\cr
\h
& $L$ && braid && $\sgn\lk(,)$ && $L\in\cQ$ &
& $L\in\cSB$ && $L\in\cB$ && $\chi_s(L)$ && $\chi_s^-(L)$ && \S\sectTab &\cr
\h
 \noalign{\hrule}\h\h
& $6_1$   && $112\1\32\3$ &&&& no && yes && yes &&  && && F &\cr\h
& $6_1^*$ &&              &&&& no && yes && yes &&  && && F &\cr\h
 \noalign{\hrule}\h\h
& $6_2$   && $111\21\2$   &&&& no && yes && yes &&      &&   && F &\cr\h
& $6_2^*$ &&              &&&& no && no  && no  && $-1$ && 1 && B &\cr\h
 \noalign{\hrule}\h\h
& $6_3$   && $11\21\2\2$  &&&& no && no  && no  && $-1$ && 1 && B &\cr\h
 \noalign{\hrule}\h\h
& L6a1(0)\;&&               &&$-$&& no && no && no && $-2$ && 2 && B &\cr\h
& L6a1(0)\*&& $\12\32\1232$ &&$+$&& no && no && no && $-2$ && 0 && B &\cr\h
& L6a1(1)\;&& $112\1211$    &&$+$&& yes&& yes&& yes&&      &&   &&   &\cr\h
& L6a1(1)\*&&               &&$-$&& no && no && no && $-2$ && 2 && B &\cr\h
 \noalign{\hrule}\h\h
& L6a2(0)\;&&               &&$-$&& no && no && no && $-2$ && 2 && B &\cr\h
& L6a2(0)\*&& $1\212211$    &&$+$&& yes&& yes&& yes&&      &&   &&   &\cr\h
 \noalign{\hrule}\h\h
& L6a3(0)\;&&               &&$-$&& no && no && no && $-4$ && 2 && B &\cr\h
& L6a3(0)\*&& $111111$      &&$+$&& yes&& yes&& yes&&      &&   &&   &\cr\h
& L6a3(1)\;&&$2\32\12\32331$&&$+$&& yes&& yes&& yes&&      &&   &&   &\cr\h
& L6a3(1)\*&&               &&$-$&& no && no && no && $ 0$ && 2 && B &\cr\h
 \noalign{\hrule}\h\h
& L6a4     && $1\21\21\2$   &&0\;0\;0&& no && yes&& yes&&    &&  && F &\cr\h
 \noalign{\hrule}\h\h
& L6a5(0,0)\;&&                &&$---$&& no && no && no && $-1$ && 3 && C &\cr\h
& L6a5(0,0)\*&&$\1\2321213\232$&&$+++$&& yes&& yes&& yes&&      &&   &&   &\cr\h
& L6a5(0,1)\;&&   $1\211\21$   &&$++-$&& no && yes&& yes&&      &&   && F &\cr\h
& L6a5(0,1)\*&&                &&$--+$&& no && no && no && $-1$ && 1 && C &\cr\h
\noalign{\hrule}\h\h
& L6n1(0,0)\;&&   $\1\12112$   &&$++-$&& yes&& yes&& yes&&      &&   &&   &\cr\h
& L6n1(0,0)\*&&                &&$--+$&& no && ?  && yes&&      &&   && [\refKO] &\cr\h
& L6n1(0,1)\;&&                &&$---$&& no && no && no && $-3$ && 3 && C &\cr\h
& L6n1(0,1)\*&&    $121212$    &&$+++$&& yes&& yes&& yes&&      &&   &&   &\cr\h
\noalign{\hrule}
}}
}
\endinsert

\midinsert
\hbox to 100mm{\hfill Table \tabB. Composite non-split links.}
\medskip
\centerline{
\vbox{\offinterlineskip
\def\h {height2pt&\omit&&\omit&&\omit&&\omit&&\omit&&\omit&&\omit&&\omit&\cr}
    
\def\1{\bar1}\def\2{\bar2}\def\3{\bar3}\def\4{\bar4}\def\*{${}^*$}

\hrule
\halign{&\vrule#&\strut\;\hfil#\hfil\,\cr
\h
& $L$ && braid && $L\in\cQ$ && $L\in\cSB$&&$L\in\cB$&&$\chi_s(L)$&&$\chi_s^-(L)$&&
                                           see \S\sectTab &\cr
\h
 \noalign{\hrule}\h\h
& $3_1\#3_1$     && $111222$ && yes && yes && yes &&      &&   &&   &\cr\h
& $3_1\#3_1^*$   &&          && no  && yes && yes &&      &&   && E &\cr\h
& $3_1^*\#3_1^*$ &&          && no  && no  && no  && $-3$ && 1 && C &\cr\h
 \noalign{\hrule}\h\h
& $4_1\#2_1$     &&  $1\21\233$  && no  && yes && yes &&    &&   && F &\cr\h
& $4_1\#2_1^*$   && $1\21\2\3\3$ && no  && no  && no  &&  0 && 2 && B &\cr\h
 \noalign{\hrule}\h\h
& L4a1(0)$\#2_1\;$   &&               && no  && no  && no  && $-1$ && 1 && B &\cr\h
& L4a1(0)\*$\#2_1$   && $2\121133$    && yes && yes && yes &&      &&   &&   &\cr\h
& L4a1(0)$\#2_1^*\;$ &&               && no  && no  && no  && $-1$ && 3 && B &\cr\h
& L4a1(0)\*$\#2_1^*$ && $2\1211\3\3$  && no  && no  && no  && $-1$ && 1 && B &\cr\h
 \noalign{\hrule}\h\h
& L4a1(1)$\#2_1\;$   &&    $111122$     && yes && yes && yes &&      &&   &&   &\cr\h
& L4a1(1)\*$\#2_1$   &&                 && no  && no  && no  && $-1$ && 1 && B &\cr\h
& L4a1(1)$\#2_1^*\;$ &&    $1111\2\2$   && no  && yes && yes &&      &&   && E &\cr\h
& L4a1(1)\*$\#2_1^*$ &&                 && no  && no  && no  && $-1$ && 3 && B &\cr\h
 \noalign{\hrule}\h\h
& $I_{+++}$ &&    $112233$    && yes && yes && yes &&      &&   &&     &\cr\h
& $I_{++-}$ &&   $1122\3\3$   && no  && yes && yes &&      &&   && A   &\cr\h
& $I_{+-+}$ &&   $11\2\233$   && no  && yes && yes &&      &&   && A,F &\cr\h
& $I_{+--}$ &&  $11\2\2\3\3$  && no  && no  && no  && $ 0$ && 2 && C   &\cr\h
& $I_{-+-}$ &&  $\1\122\3\3$  && no  && no  && no  && $ 0$ && 2 && C   &\cr\h
& $I_{---}$ && $\1\1\2\2\3\3$ && no  && no  && no  && $-2$ && 4 && B,C &\cr\h
 \noalign{\hrule}\h\h
& $Y_{+++}$ &&    $1122322\3$    && yes && yes && yes &&      &&   &&     &\cr\h
& $Y_{++-}$ &&   $11223\2\2\3$   && no  && yes && yes && $ 0$ && 2 && A   &\cr\h
& $Y_{+--}$ &&  $11\2\23\2\2\3$  && no  && no  && no  && $ 0$ && 2 && C   &\cr\h
& $Y_{---}$ && $\1\1\2\23\2\2\3$ && no  && no  && no  && $-2$ && 4 && B,C &\cr\h
\noalign{\hrule}
}}
}
\endinsert

\midinsert
\hbox to 100mm{\hfill Table \tabC. Split links.}
\medskip
\centerline{
\vbox{\offinterlineskip
\def\h {height2pt&\omit&&\omit&&\omit&&\omit&&\omit&&\omit&&\omit&&\omit&\cr}
    
\def\1{\bar1}\def\2{\bar2}\def\3{\bar3}\def\4{\bar4}\def\5{\bar5}\def\*{${}^*$}
\let\U=\sqcup

\hrule
\halign{&\vrule#&\strut\;\hfil#\hfil\,\cr
\h
& $L$ && braid && $L\in\cQ$ && $L\in\cSB$&&$L\in\cB$&&$\chi_s(L)$&&$\chi_s^-(L)$&&
                                           see \S\sectTab &\cr
\h
 \noalign{\hrule}\h\h
& $3_1\U 3_1$     && $111\,333$ && yes && yes && yes &&      &&   &&  &\cr\h
& $3_1\U 3_1^*$   &&            && no  &&  ?  && yes &&      &&   && E &\cr\h
& $3_1^*\U 3_1^*$ &&            && no  && no  && no  && $-2$ && 2 && C &\cr\h
 \noalign{\hrule}\h\h
& $4_1\U 2_1$     &&  $1\21\2\,44$  && no  && no  && no  && $-1$ && 1 && D &\cr\h
& $4_1\U 2_1^*$   && $1\21\2\,\4\4$ && no  && no  && no  && $-1$ && 3 && D &\cr\h
 \noalign{\hrule}\h\h
& L4a1(0)${}\U 2_1\;$   &&                 && no  && no  && no  &&  0 && 2 && D &\cr\h
& L4a1(0)\*${}\U 2_1$   && $2\1211\,44$    && yes && yes && yes &&    &&   &&   &\cr\h
& L4a1(0)${}\U 2_1^*\;$ &&                 && no  && no  && no  &&  0 && 4 && C &\cr\h
& L4a1(0)\*${}\U 2_1^*$ && $2\1211\,\4\4$  && no  && no  && no  &&  0 && 2 && D &\cr\h
 \noalign{\hrule}\h\h
& L4a1(1)${}\U 2_1\;$   &&  $1111\,33$     && yes && yes && yes &&      &&   &&   &\cr\h
& L4a1(1)\*${}\U 2_1$   &&                 && no  && no  && no  && $-2$ && 0 && D &\cr\h
& L4a1(1)${}\U 2_1^*\;$ &&  $1111\,\3\3$   && no  && no  && yes && $-2$ && 0 && E &\cr\h
& L4a1(1)\*${}\U 2_1^*$ &&                 && no  && no  && no  && $-2$ && 2 && C &\cr\h
 \noalign{\hrule}\h\h
& $2_1\#2_1\U 2_1\;$      &&    $1122\,44$    && yes&& yes&& yes&&      &&   &&   &\cr\h
& $2_1\#2_1^*\U 2_1\;$    &&   $11\2\2\,44$   && no && yes&& yes&&      &&   && A &\cr\h
& $2_1^*\#2_1^*\U 2_1\;$  &&  $\1\1\2\2\,44$  && no && no && no && $-1$ && 3 && D &\cr\h
& $2_1\#2_1\U 2_1^*\;$    &&   $1122\,\4\4$   && no && no && yes&& $-1$ && 1 && A &\cr\h
& $2_1\#2_1^*\U 2_1^*\;$  &&  $11\2\2\,\4\4$  && no && no && no && $ 1$ && 3 && C &\cr\h
& $2_1^*\#2_1^*\U 2_1^*\;$&& $\1\1\2\2\,\4\4$ && no && no && no && $-1$ && 5 && C &\cr\h
 \noalign{\hrule}\h\h
& $2_1\U 2_1\U 2_1\;$      &&   $11\,33\,55$   && yes && yes && yes&&   &&   &&   &\cr\h
& $2_1\U 2_1\U 2_1^*\;$    &&  $11\,33\,\5\5$  && no  && no  && yes&& 0 && 2 && A &\cr\h
& $2_1\U 2_1^*\U 2_1^*\;$  && $11\,\3\3\,\5\5$ && no  && no  && no && 0 && 4 && D &\cr\h
& $2_1^*\U 2_1^*\U 2_1^*\;$&&$\1\1\,\3\3\,\5\5$&& no  && no  && no && 0 && 6 && C &\cr\h
\noalign{\hrule}
}}
}
\endinsert


\subhead\sectIY. The links $2_1\# 2_1\# 2_1$, $2_1\# 2_1\# 2_1^*$, etc
\endsubhead
Among the links with 6 crossings, almost all connected sums are uniquely determined by
their summands. The only exceptions are
connected sums of three copies of the Hopf link $2_1$ and/or its mirror image $2_1^*$.
Such a sum is determined by the following graph with signed edges.
Its vertices correspond to the link components. Two vertices are connected with
an edge if the corresponding components are linked,
and this edge is labeled by the linking number,
which is always $\pm1$. This graph is a tree.
There are only two 4-vertex trees: an $I$-shaped tree (a 4-path) and a $Y$-shaped
tree (with a vertex of degree 3). We denote the corresponding links by
$I_{\pm\pm\pm}$ and $Y_{\pm\pm\pm}$ respectively. The signs in the notation
$I_{\pm\pm\pm}$ appear in the same order as in the graph. 

Any connected sum of two strong $\C$-boundaries is a strong $\C$-boundary,
and $2_1\# 2_1^*$ is a strong $\C$-boundary (see Example~\exaQQ). Hence
all links 
of the form $2_1\# 2_1\# 2_1^*$ (i.e., $I_{++-}$, $I_{+-+}$, and $Y_{++-}$)
are strong $\C$-boundaries. Note that $I_{++-}$ and $Y_{++-}$ are strong $\C$-boundaries
by [\refKO, Theorem~2.1], and we give another
strong $\C$-boundary realization of $I_{+-+}$ in Figure~\figMain.


\subhead\sectZapret. Restrictions on $\C$-boundaries
\endsubhead
Our proofs that some links are not (strong) $\C$-boundaries are based on
Kronheimer-Mrowka's Theorem [\refKM]
(former Thom Conjecture) and its version for
immersed 2-surfaces in $\CP^2$ with negative double points
(see [\refFS], [\refM, \S2]), which was actually proven in [\refKM]
but was not explicitly formulated.
 
\bigskip
\proclaim{ Theorem \thKM } {\rm(The Immersed Thom Conjecture.)}
Let $\Sigma$ be a connected oriented closed surface of genus $g$ and
$j:\Sigma\to\CP^2$
be an immersion which has only negative ordinary double points as self-crossings.
Let $j_*([\Sigma])=d[\CP^1]\in H_2(\CP^2)$ with $d>0$. Then $g$ is bounded below by the
genus of a smooth algebraic curve of degree $d$, that is,
$g \ge (d-1)(d-2)/2$.
\endproclaim

Given a link $L$ in $S^3=\partial B^4$, we define the {\it slice
Euler characteristic} of $L$ by $\chi_s(L)=\max_\Sigma\chi(\Sigma)$ where the
maximum is taken over all embedded smooth oriented surfaces $\Sigma$ without
closed components, such that $\partial\Sigma=L$.
Similarly, we define the {\it slice negatively immersed
Euler characteristic} of $L$ by $\chi_s^-(L)=\max_{(\Sigma,j)}\chi(\Sigma)$ where the
maximum is taken over all immersions $j:(\Sigma,\partial\Sigma)\to(B^4,S^3)$
of oriented surfaces $\Sigma$ without closed components such that
$j(\Sigma)$ has only negative double points and $j(\partial\Sigma)=L$.
Theorem \thKM\ easily implies the following two statements.

\proclaim{ Proposition \propB } {\rm(See [\refKO, Prop.~3.2 and Lem.~3.8].)}
If $\chi_s(L) < \chi_s^-(L)$, then $L\not\in\cSB$. If, moreover,
$L$ does not have a proper sublink which has zero linking number with its
complement (for example, if $L$ is a knot),
then $L\not\in\cB$.
\endproclaim

We see in [\refKO, Table 1] and in Tables \tabA--\tabC\ here that for any link with
$\le 6$ crossings with a single exception $3_1^*\#2_1$, we have
either $L\in\cSB$ or $\chi_s(L)<\chi_s^-(L)$.

\proclaim{ Proposition \propC } {\rm(See [\refKO, Prop.~3.7]; cf.~[\refBR, Cor.~1.6].)}
If $L$ is a strong $\C$-boundary and $-L^*$ is a
(not necessarily strong) $\C$-boundary, then $\chi_s(L)=\chi_s^-(L)\ge 1$.
\endproclaim

Now let us treat the cases which are not covered by Propositions~\propB, \propC.

\proclaim{ Proposition \propD } The links marked by ``D'' in the last column
of Table~\tabC\ are not $\C$-boundaries.
\endproclaim

\demo{ Proof }
{\it Case 1. The links $4_1\sqcup 2_1$ and $4_1\sqcup 2_1^*$.}
Suppose that there exists a $\C$-boundary realization $L=A\cap B$ of one of these links.
Let $L=L_4\cup L_2$ where $L_4$ is $4_1$ and $L_2$ is $2_1$ or $2_1^*$.
$L$ is not a strong $\C$-boundary by Proposition~\propB. Hence
$A\setminus B$ is a disjoint union $A_2\cup A_4$ with $\partial A_k=-L_k$,
where either $A_4$ or $A_2$ is bounded.
The interior piece $A_0=A\cap B$ of $A$ is connected because otherwise
$A$ would be disconnected (here we use the fact that there is only one partition
of $L$ into sublinks with zero linking number). 
Hence $\chi(A_0)\le-1$. Let us replace $A_0$ by $A'_2\sqcup A'_4$, $\partial A'_k=L_k$,
where $A'_2$ is an embedded annulus and $A'_4$ is an immersed disk with a negative
self-crossing. Recall that $A_2$ or $A_4$ is bounded. If it is $A_2$,
then $A_4\cup A'_4$ contradicts Theorem~\thKM; if it is $A_4$, then so does
$A_2\cup A'_2$ because in this case
$\chi(A_4)\le-1$ ($4_1$ is not slice), and hence the replacement $A_0\cup A_4\to A'_2$
decreases $\chi_s^-$ (see Figure~\figD).

\midinsert
\centerline{\epsfxsize=120mm\epsfbox{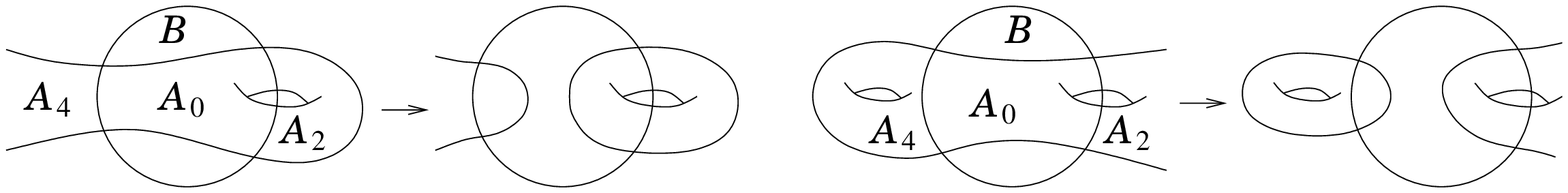}}
\botcaption{Figure \figD} Proof that $4_1\sqcup 2_1$ and $4_1\sqcup 2_1^*$
 are not $\C$-boundaries.
\endcaption
\endinsert

\smallskip
{\it Case 2. The links} L4a1(0)${}\sqcup 2_1$, L4a1(0)${}^*\sqcup 2_1^*$,
 {\it and} L4a1(1)${}^*\sqcup 2_1$.
The same proof as in the previous case; see Figure~\figDD.

\midinsert
\centerline{\epsfxsize=55mm\epsfbox{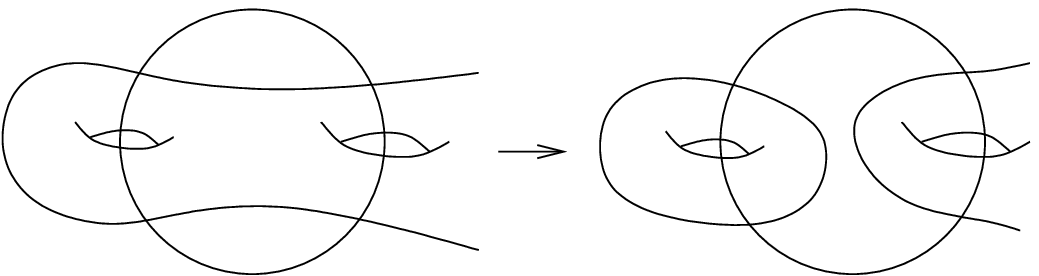}}
\botcaption{Figure \figDD} Proof that L4a1(0)${}\sqcup 2_1$, L4a1(0)${}^*\sqcup 2_1^*$,
  L4a1(1)${}^*\sqcup 2_1\not\in\cB$.
\endcaption
\endinsert

\medskip
{\it Case 3. The link $L=2_1^*\# 2_1^*\sqcup 2_1$.}
Suppose that $L=A\cap B$ is a $\C$-boundary realization of $L$.
We have $\chi_s(L)=-1$ and $\chi_s^-(L)=3$. Thus, if we replace $A\cap B$
with a negatively self-crossing surface $A'_0$ which realizes $\chi_s^-(L)=3$,
we obtain a contradiction with Theorem~\thKM\ unless at least two spheres split out.

Note that if we glue together several surfaces along their boundaries, we
may obtain a sphere only if at least two of the glued surfaces are disks.
The link $L$ does not have any component which has zero linking number with its
complement, hence no component of $A\setminus B$ is a disk.
The negatively self-crossing surface $A'_0$ has at most three disks,
hence at most one sphere can split out.

\medskip
{\it Case 4. The link $L=2_1\sqcup 2_1^*\sqcup 2_1^*$.}
Suppose that $L=A\cap B$ is a $\C$-boundary realization of $L$.
We have $\chi_s(L)=0$ and $\chi_s^-(L)=4$. Thus, as in Case 3,
if we replace $A\cap B$ with a negatively self-crossing surface $A'_0$,
at least two spheres must split out. Again as in Case 3, no component of $A\setminus B$
is a disk but this time $A'_0$ may have four disks (bounded by the components
of $2_1^*\sqcup 2_1^*$).
Thus all these disks are present in $A'_0$, all of them
make part of the two spheres which split out, and $\chi(A\cap B)=0$
(if $\chi(A\cap B)$ were negative, then one more sphere should split out).
This may happen only when
$A\setminus B=A_1^-\cup A_2^-\cup A^+$ where $A_j^-$ are cylinders, $A^+$ is the
unbounded component, each $\partial A_j^-$ is $2_1^*$, and $\partial A^+$ is $2_1$.

The condition $\chi(A\cap B)=0$ implies that $A\cap B$ is a union of three cylinders.
Since all the linking numbers between their boundaries are zero, one of them must
be bounded by $2_1^*$ (while a priori the boundaries of the other two may realize any
$(2,2)$-partition of the four components of $2_1\sqcup 2_1^*$).
This fact contradicts the connectedness of $A$.
\qed\enddemo

\midinsert
\centerline{\lower-6mm\hbox to15mm{$\,\,\,6_1$}\epsfxsize=90mm \epsfbox{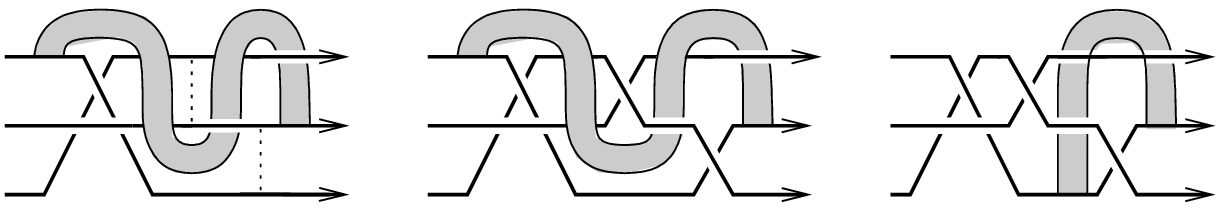}}\bigskip
\centerline{\lower-5mm\hbox to15mm{$\,\,\,6_1^*$}\epsfxsize=90mm\epsfbox{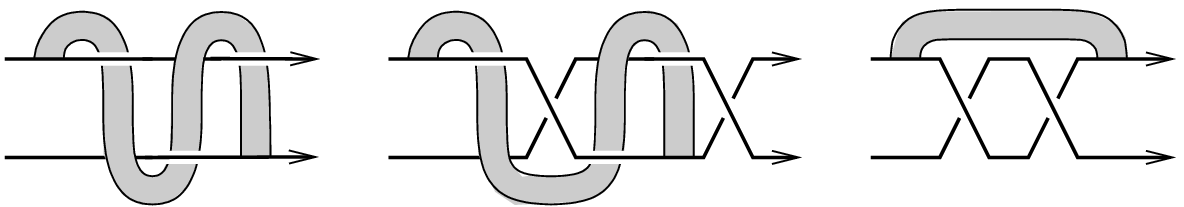}}\bigskip
\centerline{\lower-4mm\hbox to15mm{$\,\,\,6_2$}\epsfxsize=90mm \epsfbox{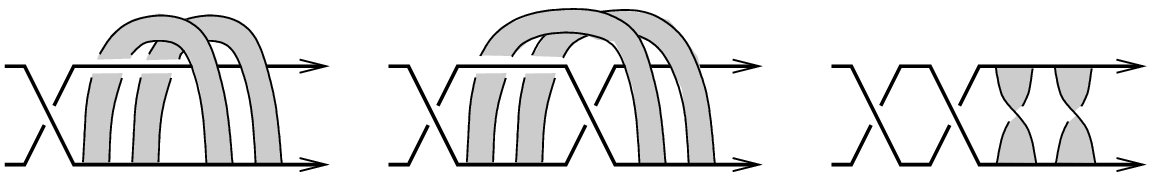}}\bigskip
\centerline{\lower-6mm\hbox to15mm{$4_1\#2_1$}\epsfxsize=90mm \epsfbox{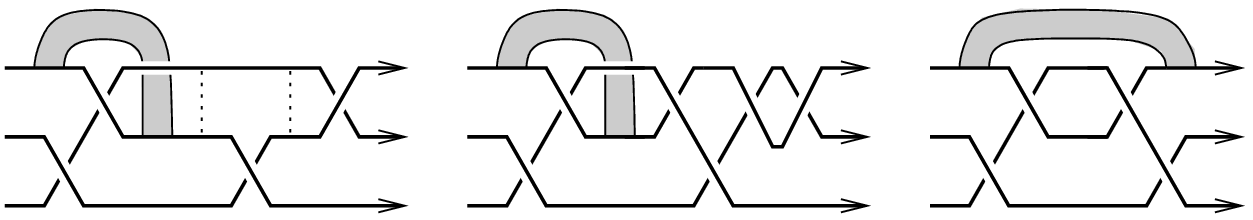}}\bigskip
\centerline{\lower-7mm\hbox to15mm{L6a4}\epsfxsize=90mm \epsfbox{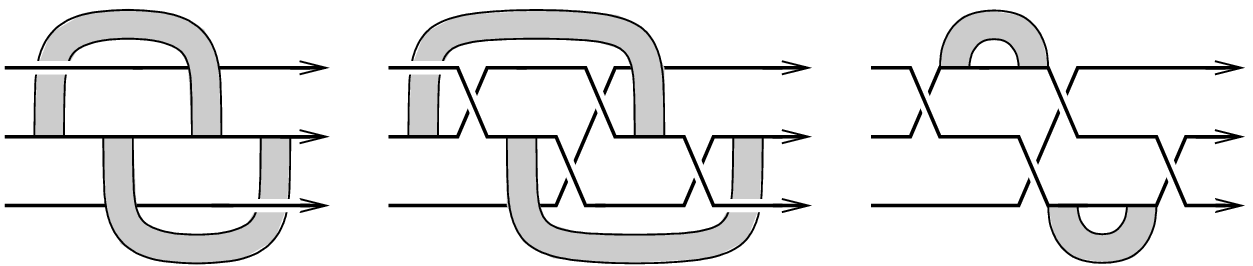}}\bigskip
\centerline{\hskip-3mm\lower-5mm\hbox to17mm{L6a5(0,1)}\epsfxsize=90mm\epsfbox{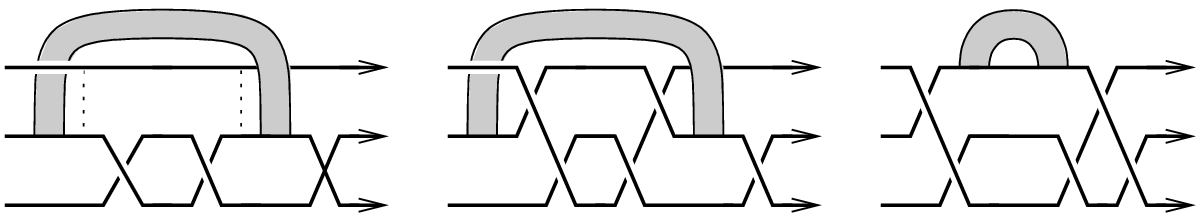}}
\bigskip\centerline{\lower-6mm\hbox to15mm{$I_{+-+}$}\epsfxsize=90mm \epsfbox{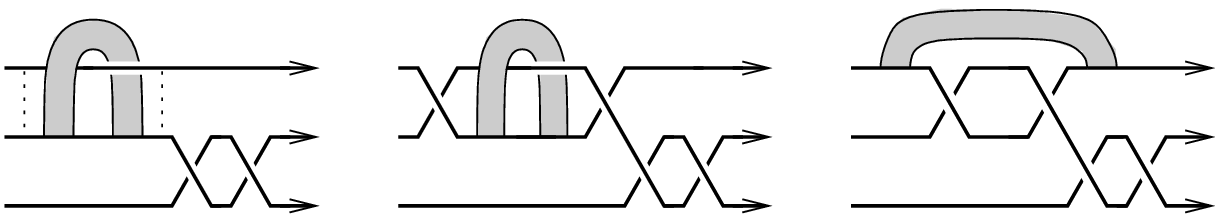}}
\botcaption{Figure \figMain}
A strong $\C$-boundary realization of the knots $6_1$, $6_1^*$, $6_2$, $4_1\# 2_1$
  and the links L6a4 (Borromeo rings), L6a5(0,1), and $I_{+-+}$.
\endcaption
\endinsert


\subhead\sectTab. Comments on Tables~\tabA--\tabC
\endsubhead
The tables are organized similarly to Table 1 in [\refKO].
The links are named according to [\refLMknot], [\refLMlink]
except the two-component Hopf links $2_1$, $2_1^*$, and their connected sums, whose
notation is explained in \S\sectIY.
The braid notation and the signs
of linking numbers help 
to identify the links.
The invariants $\chi_s(L)$ and $\chi_s^-(L)$ are introduced in \S\sectZapret.
They are computed by standard methods (some of them are mentioned in [\refKO, \S6])
and we omit the details.
The same for the quasipositivity: in all the considered cases,
if $L\in\cQ$, this fact is clear from the braid given in the tables, and if $L\not\in\cQ$,
this fact follows either from [\refO], or from the Franks-Williams-Morton Inequality
[\refFW], [\refMo] in the form given in [\refBR, Thm.~3.2] (see [\refKO, Thm.~6.1]).
In all the cases when $L\in\cB\setminus\cSB$, the fact that $L\not\in\cSB$ follows
from the inequality $\chi_s(L)<\chi_s^-(L)$ (see Proposition~\propB).

The letters A--F in the last column of each table mean the following:
\roster
\item"A." $L\in\cB$ or $L\in\cSB$ by the additivity of (strong)
         $\C$-boundaries under split or connected sums [\refKO, Prop.~3.6]
         (see also \S\sectIY).
\item"B." $L\not\in\cB$ by Proposition \propB.
\item"C." $L\not\in\cB$ by Proposition \propC.
\item"D." $L\not\in\cB$ by Proposition \propD.
\item"E." See Examples~\exaQQ\ and \exaQQii.
\item"F." See Figure~\figMain.
\endroster


\head\sectConclude. Concluding remarks
\endhead


\subhead\sectSqueezed. Squeezed knots
\endsubhead
The referee drew my attention to a very interesting paper [\refFLL] by
Feller, Lewark, and Lobb, where they introduce 
{\it squeezed knots}. These are knots which appear as a slice
of a genus-minimizing oriented connected smooth cobordism between a positive
torus knot and a negative torus knot.

\proclaim{ Proposition \propSqueezed } Any $\C$-boundary knot is squeezed.
\endproclaim

\demo{ Proof } Let $A\cap B$, $A=\{f(z,w)=0\}$, be a $\C$-boundary realization of a knot $K$. Let $n=\deg f$.
Without loss of generality we may assume that $A$ is smooth, $f(0,0)=0$, and $\deg_z f=n$.
Let $A_\varepsilon = \{f(z,w)=\varepsilon w^{n+1}\}$. If $|\varepsilon|$ is sufficiently small,
then $K_\varepsilon=A_\varepsilon\cap B$ is isotopic to $K$ and
$A_\varepsilon\cap(B_R\setminus B_r)$, $0<r\ll 1\ll R$,
is a cobordism between the positive torus knot $T(n,n+1)$ and the unknot. This cobordism is
genus-minimizing by Kronheimer-Mrowka's theorem.
\qed\enddemo

As proven in [\refFLL], the knots $9_{42}$, $10_{125}$, $10_{132}$, and $10_{136}$ are not squeezed.
Hence they and their mirrors are not $\C$-boundaries.
For these knots themselves, this fact follows from Proposition~\propB, because they are not slice
(see [\refLMknot]) but $\chi_s^-=1$. The latter fact can be checked using the braids
(see the beginning of \S6 in [\refKO]):
$$
\split
   \chi_s^-(9_{42})&=\chi_s^-(\bar1\bar1\bar1211\bar32\bar3)
                  \ge\chi_s(\bar1(\bar1\bar1211)(\bar323))=1,
\\
   \chi_s^-(10_{125})&=\chi_s^-(\bar1\bar1\bar1\bar1\bar121112)
                    \ge\chi_s((\bar1\bar1\bar12111)2)=1,
\\
   \chi_s^-(10_{132})&=\chi_s^-(111\bar2\bar1\bar1\bar2\bar32\bar3\bar3)
                       \ge\chi_s(1(11\bar2\bar1\bar1)(\bar2\bar32))=1,
\\
   \chi_s^-(10_{136})&=\chi_s^-(\bar1 \bar1 \bar2 3\bar2 1 \bar2 \bar2 3 2 2)
                        \ge\chi_s((2 3\bar2) 1 (\bar2 \bar2 3 2 2))=1.
\endsplit
$$

In contrary, the methods of our paper are not sufficient to prove that the mirror images of these four knots
are not $\C$-boundaries. In the first arxiv version of this paper I asked whether $\chi_s^-(K)=\chi_s(K)$
implies that $K$ is $\C$-boundary. The knots $9_{42}^*$, $10_{125}*$, $10_{136}^*$ provide a negative answer.
Indeed, if $K$ is one of them, then the signature of $K\# 2_1$ equals $3$ which implies that
$\chi_s^-(K)<0$, hence $\chi_s^-(K)=\chi_s(K)=-1$.

Notice that a negative answer to the same question for links is provided by the link $L=3_1^*\sqcup 2_1$.
We have $\chi_s(L)=\chi_s^-(L)=0$ whereas $L\not\in\cB$ by Proposition~\propC\ (see [\refKO]).
It also seems plausible that $\chi_s^-(K)=\chi_s(K)=-1$ for $K=8_{18}$ (the braid closure of
$(\sigma_1\sigma_2^{-1})^4$), whereas $K=K^*$ and hence $K$ is not a $\C$-boundary by [\refBR, Cor.~1.6].

\if01{
(* 10_130 *)
f[4,{1,1,1,-2,-1,-1,-2,-2,-3,2,-3}]         (* 0,0 *)
f[5,{1,1,1,-2,-1,-1,-2,-2,-3,2,-3, 4,4}]    (*-1,0 *)
f[5,{1,1,1,-2,-1,-1,-2,-2,-3,2,-3, -4,-4}]  (* 1,0 *)

(* 9_42 *)
f[4,{-1,-1,-1,2,1,1,-3,2,-3}]             (* -2,0 *)
f[5,{-1,-1,-1,2,1,1,-3,2,-3,4,4}]         (* -1,0 *)
f[5,{-1,-1,-1,2,1,1,-3,2,-3,-4,-4}]       (* -3,0 *)

(* 10_125 *)
f[3,{-1,-1,-1,-1,-1,2,1,1,1,2}]             (* -2,0 *)
f[4,{-1,-1,-1,-1,-1,2,1,1,1,2,3,3}]         (* -1,0 *)
f[4,{-1,-1,-1,-1,-1,2,1,1,1,2,-3,-3}]       (* -3,0 *)

(* 10_132 *)
f[4,{1,1,1,-2,-1,-1,-2,-3,2,-3,-3}]         (* 0,0 *)
f[5,{1,1,1,-2,-1,-1,-2,-3,2,-3,-3,4,4}]     (*-1,0 *)
f[5,{1,1,1,-2,-1,-1,-2,-3,2,-3,-3,-4,-4}]   (* 1,0 *)

(* 10_136 *)
f[4,{-1,-1,-2,3,-2,1,-2,-2,3,2,2}]             (* -2,0 *)
f[5,{-1,-1,-2,3,-2,1,-2,-2,3,2,2,4,4}]         (* -1,0 *)
f[5,{-1,-1,-2,3,-2,1,-2,-2,3,2,2,-4,-4}]       (* -3,0 *)

}\fi


\subhead\sectSliceKnots. Slice knots
\endsubhead
Up to now, all known proofs that a certain knot is not $\C$-boundary,
automatically gives that all concordant knot are not $\C$-boundaries neither.
Thus there are no tools to prove that a slice knot is not a $\C$-boundary.

There are 15 slice knots with $\le9$ crossings:
$6_1$, $8_8$, $8_9$, $8_{20}$, $9_{27}$, $9_{41}$, $9_{46}$, their
mirror images ($8_9$ is amphicheiral), $3_1\#\, 3_1^*$, and $4_1\#\,4_1$.
Among them, the following ten are realized as $\C$-boundaries:
$6_1$, $6_1^*$ (Figure~\figMain);
$8_8$, $8_9$, $9_{27}^*$, $9_{41}^*$ (Figure~\figSliceKnots);
$8_{20}^*$, $9_{46}^*$ (quasipositive knots); $3_1\#\,3_1^*$ and
$4_1\#\,4_1$ (Examples~\exaQQ\ and \exaQQiii). The realizability of
$8_8^*$, $8_{20}$, $9_{27}$, $9_{41}$, and $9_{46}$ is unknown.

\midinsert
\centerline{\lower-6mm\hbox to10mm{$8_8$}\epsfxsize=98mm\epsfbox{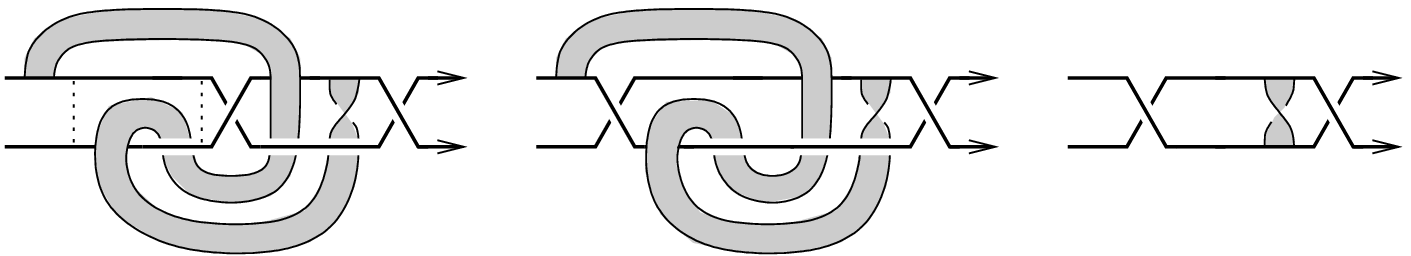}}\bigskip
\centerline{\lower-4mm\hbox to10mm{$8_9$}\epsfxsize=98mm\epsfbox{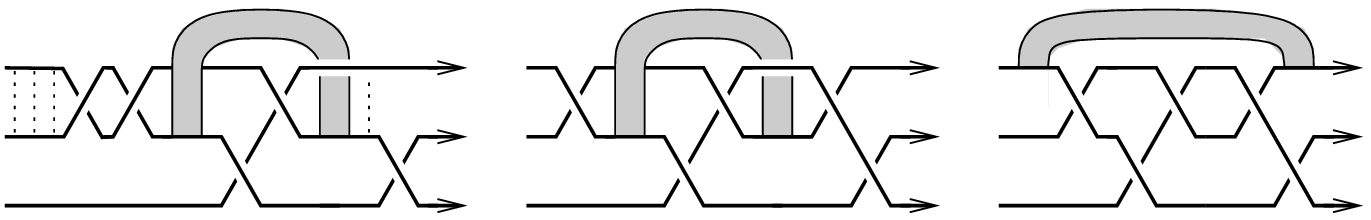}}\bigskip
\centerline{\lower-7mm\hbox to10mm{$9_{27}^*$}\epsfxsize=98mm \epsfbox{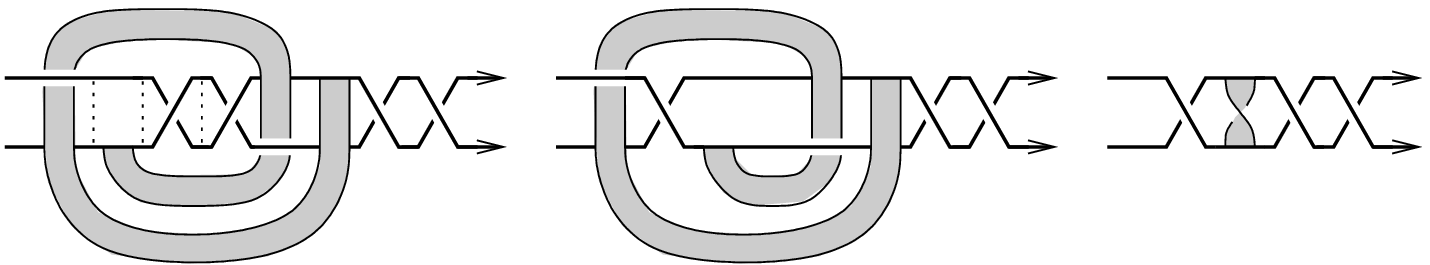}}\bigskip
\centerline{\lower-8mm\hbox to6mm{$9_{41}^*$}\epsfxsize=112mm\epsfbox{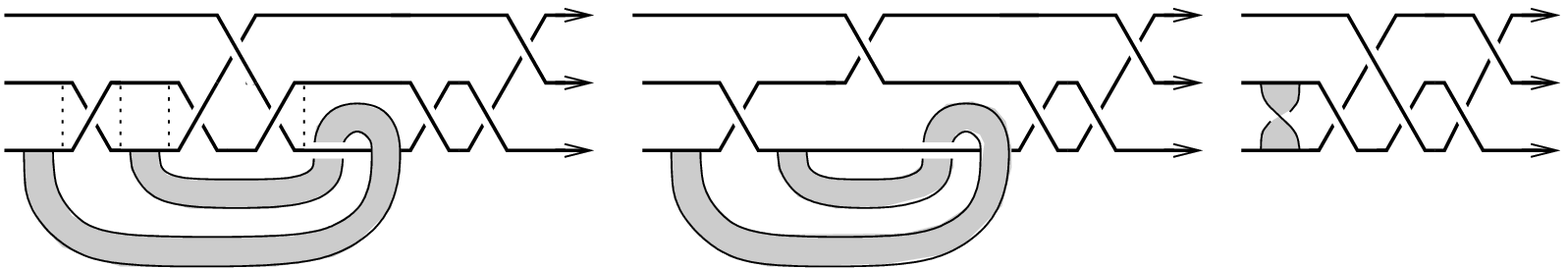}}
\botcaption{Figure \figSliceKnots}
A $\C$-boundary realization of the knots $8_8$, $8_9=8_9^*$,
     $9_{27}^*$, $9_{41}^*$.
\endcaption
\endinsert


\subhead\sectQuestions. Some open questions
\endsubhead
The study of $\C$-boundary links is at the very beginning. So, most
questions on this subject are open. However, I would like to emphasize some of them.

\medskip\noindent
{\bf Question \queConstr.} Can any $\C$-boundary be realized by the
construction from \S\sectConstr?

\if01{
\medskip\noindent
{\bf Question \queChi.} Does $\chi_s(K)=\chi_s^-(K)$ imply
that a knot $K$ is $\C$-boundary?

\medskip
A candidate for a negative answer is the knot $K=8_{18}$
(the braid closure of the 3-braid $(\sigma_1\sigma_2^{-1})^4$).
Since $K=-K^*$, it 
is not a $\C$-boundary by [\refBR, Cor.~1.6]. One has $\chi_s(K)=-1$
and it seems plausible that $\chi_s^-(K)=-1$.

The answer to this question for links is negative: we have
$\chi_s(L)=\chi_s^-(L)=0$ for $L=3_1^*\#2_1$ whereas $L\not\in\cB$ by
Proposition~\propC\ (see [\refKO]).
}\fi


\medskip\noindent
{\bf Question \queSliceI.} Do there exist slice knots which are not $\C$-boundaries?
In particular, does there exist a knot $K$ such that $K\,\#\,{-K^*}$
is not a $\C$-boundary?

\medskip
As we pointed out in \S\sectSliceKnots, if the answer is affirmative, then essentially
new ideas should be evoked to prove this fact.

\medskip\noindent
{\bf Question \queSliceII.} Is it true that for any slice knot $K$,
either $K$ or $K^*$ is $\C$-boundary?

\medskip
The answer is affirmative for
knots with at most 9 crossings (see \S\sectSliceKnots).

\medskip\noindent
{\bf Question \queStrong.} Is it true that $K\sqcup -K^*$ is not a strong $\C$-boundary
for any non-trivial quasipositive knot $K$?

\medskip\noindent
{\bf Question \queCancel.} (Cancellation of unknot in a split sum.)
If $L\sqcup O$ is a (strong) $\C$-boundary, is $L$ a (strong) $\C$-boundary?

\medskip\noindent
{\bf Question \queSplit.} Do there exist non-split non-strong $\C$-boundaries?

\medskip
The link L6n1(0,0)${}^*$ (see [\refKO, \S4]) is a candidate for
an affirmative answer.


\enddocument